\documentclass[a4paper, 12pt, authoryear]{elsarticle}
\journal{arXiv}
\usepackage[unicode=true]{hyperref}
\usepackage{cite}
\usepackage{amsthm}
\usepackage{graphics,graphicx,epsfig}
\usepackage{cite}
\usepackage{calc}
\usepackage{hyperref}
\usepackage{latexsym}
\usepackage{amssymb,amstext,amsmath,amsfonts}

\usepackage{orcidlink}
\usepackage{color}
\usepackage{soul}
\usepackage[margin=1.25in]{geometry}

\def\Re{\mathbb{R}}
\DeclareMathOperator{\im}{im}
\DeclareMathOperator{\sign}{sign}

\DeclareMathOperator{\de}{\mathbf{d}}

\newcommand{\norm}[1]{\left\lVert#1\right\rVert}

\theoremstyle{plain}
\newtheorem{theorem}{Theorem}[section]
\newtheorem{lemma}{Lemma}[section]
\newtheorem{proposition}{Proposition}[section]

\theoremstyle{definition}
\newtheorem{definition}{Definition}[section]

\theoremstyle{remark}

\begin{document}
\hyphenation{build}


\begin{frontmatter}

\title{A discrete wedge product on general polygonal meshes}

\author{Lenka Pt\'a\v ckov\'a}
\ead{lenka@kam.mff.cuni.cz}
\affiliation{
organization={Department of Numerical Mathematics, Faculty of Mathematics and Physics, Charles University},
addressline={Sokolovská 49/83},
city={Prague 8},
postcode={186 75},
country={Czech Republic}
}

\begin{abstract}
Discrete exterior calculus offers a coordinate--free discretization of exterior calculus especially suited for computations on meshes over curved manifolds. The discretization of the wedge product, that would be compatible with discrete exterior derivative, has been a challenging task. The cup product of cochains is traditionally considered to be an appropriate discrete wedge product. However, only the case of pure triangle or pure quadrilateral surface meshes has been studied thoroughly. In this work, we extend this tradition to general polygonal meshes. Specifically, we present explicit formulas for calculation of a cup/discrete wedge product on surface meshes that correspond to 2--dimensional pseudomanifolds, whose 2--dimensional faces are any simple polygons. We rigorously prove that the proposed product satisfies the definition of an abstract cup product; notably, we show that the product is compatible with the discrete exterior derivative in the sense that it satisfies the Leibniz product rule. Furthermore, the product is associative on the cohomology level, but not on the cochain level in general. We analyze the lack of associativity on the cochain level and prove that the error tends to zero under refinement of the mesh. We thus argue that the proposed product is an appropriate discretization of the wedge product of differential forms on general polygonal meshes.
\end{abstract}


\end{frontmatter}

\section{Introduction}\label{sec:Intro}

The exterior calculus of differential forms, introduced by \'{E}lie Joseph Cartan (1869 -- 1951), has become the foundation of modern differential geometry. It is a coordinate--free calculus and therefore greatly simplifies the analysis and calculations on curved spaces of differential manifolds.

Discrete versions of exterior calculus on surface meshes have gained increased attention in the last few decades, see \citet{Hirani2003, FEEC, Desbrun2008, PtackovaVelho2021} and references therein. Discete exterior calculus (DEC) and finite element exterior calculus (FEEC) have found applications in computational fluid dynamics \citep{Mohamed2016, Nitschke2017}, computational electromagnetism \citep{Zhang2021, Monkola2023, Raty2025}, computer graphics \citep{Desbrun2003, PtackovaVelho2021, FengGillespieCrane2023}, and many other areas in physics and geometry \citep{ArnoldHu2021}. Specifically, the wedge/discrete cup product has been employed for formulations of transport phenomena on cell complexes \citep{Berbatov2026}, for computing cohomology rings of 3D digital images \citep{GonzalezDiaz2011a}, and in definition of Lie derivative and Lie advection in \citet{PtackovaVelho2021}. For further applications, see the references in the cited articles.

On the other hand, there is an increased interest in discrete differential geometry on general polygonal meshes, see for example \citet{DeGoes2020, PtackovaVelho2021, Bunge2024, PtackovaOutrata2026} and references therein.

In this paper, we prove that the discrete wedge product $\wedge$ of discrete differential forms (cochains) on general polygonal meshes, presented in \citet{PtackovaVelho2021}, satisfies the Leibniz product rule with discrete exterior derivative $\de$ (coboundary operator). Furthermore, we show that the product corresponds to a cup product of cochains on orientable 2--dimensional pseudomanifolds, whose highest dimensional cells are any simple polygons. Even though several discretizations of wedge product on simplicial meshes have been studied previously, to the best of our knowledge, similar results for product on general polygonal meshes have not been published yet. We thus partially extend the findings of \citet{Schubel2024} from simplicial to general polygonal setting. Furthermore, we extend the result of almost associativity of the cup product of simplicial cochains \citet[Theorem 5.9]{Wilson2007} to polygonal cochains.

On smooth manifolds, the wedge product allows for building higher degree forms from lower degree ones. Similarly a cup product is a product of two cochains of arbitrary degree $p$ and $q$ that returns a cochain of degree $p+q$ located on $(p+q)$--dimensional cells. The cup product was introduced by J. W. Alexander, E. \v Cech, and H. Whitney~\citep{Whitney1957} in 1930's and it became a well--studied notion in algebraic topology, mainly in the simplicial setting \citep{Munkres1984, Fenn1983}. Later, the cup product was extended also to $n$--cubes~\citep{Massey1991}. In \citep{Arnold2012} various cup products on simplicial and cubical $n$--pseudomanifolds are presented. \citet{Arnold2012} also defines L${}^2$ cubical Whitney forms and a cubical cup product that fit together with the wedge product of cubical Whitney forms.

In \citep{GonzalezDiaz2011a} cup products on cubical complexes have been applied to compute cohomology rings of 3D digital images (represented as cubical complexes $Q(I)$). This information then allows them to simplify the combinatorial structure of $Q(I)$ and obtain a homeomorphic cellular complex $P(I)$ with fewer cells, which improves the computational efficiency.  \citet{GonzalezDiaz2011a} think of the cup product as the way the $n$--dimensional holes obtained in homology are related to each other. It is known that two objects with non isomorphic cohomology rings are not homotopic.

Formulas for a diagonal approximation on a general polygon that are used to compute cup products on the cohomology of a polygonal cellular complexes are presented in \citep{GonzalezDiaz2011b}. Their cup product uses a complex structure called the AT--model \citep[Theorem 2]{GonzalezDiaz2011b}, and the coefficient of their cochains lie in the field $\mathbb{Z}_2$. Whereas the coefficients of our cochains lie in $\Re$.

In common to previous approaches, this discretization of the wedge product is metric--independent and satisfies the core properties such as the Leibniz product rule, skew--commutativity, and associativity on closed forms.


\subsection*{Results}
We present a new cup product on polygonal 2--pseudomanifolds (Definition \ref{def:PolygonalCup}) and show that several previously studied cup products on simplicial and cubical 2--pseudomanifolds are special cases of this more general polygonal cup product. We proceed by proving that the polygonal cup product satisfies the Leibniz product rule with discrete exterior derivative (Lemma \ref{lemma:LeibnizPolygonalCup}) and it satisfies defining properties of an abstract cup product on cellular complexes (Theorem \ref{theo:PolygonalCupIsCup}). Moreover, in Proposition \ref{prop:PolygonalCup} we show that the polygonal cup product shares several further properties with the wedge product of differential forms, such skew--commutativity. In the Proposition \ref{prop:PolygonalCup} we also state under what conditions the product is associative and we prove that the error caused by the general lack of associativity tends to zero under mesh refinement.

In Section \ref{sec:discreteSetting} we review the appropriate discrete setting, concretely we define orientable pseudomanifolds, which is the domain of our interest, and we revise the boundary operator of chains of cells. In Section \ref{sec:DiscreteForms} we formally introduce discrete differential forms as cochains and discrete exterior derivative as coboundary operator. Section \ref{sec:CupProduct} then contains our the actual results announced in the previous paragraph.


\section{Discrete setting}\label{sec:discreteSetting}

We work with a general polygonal mesh, i.e., a surface mesh whose faces are arbitrary simple polygons, possibly non-convex and non-planar. From the perspective of algebraic topology such a mesh corresponds to a 2--dimensional CW complex with some additional structures: a 2--dimensional pseudomanifold. In this section we briefly revise several definitions that ultimately leads to the definition of an orientable pseudomanifold and the boundary operator.

A CW complex is a cellular complex, invented by J. H. C. Whitehead, and it is a more general complex than a simplicial or a singular complex. According to \citep{Munkres1984}, cellular complexes are more effective for computations of homology groups of given spaces. The following set of definitions were adapted from definitions in \citep[Paragraph 38]{Munkres1984}.

\begin{definition}
 A space is called a \textbf{\textit{m}-cell} if it is homeomorphic to the unit $m$-dimensional $B^m$. It is called an \textbf{open cell} of dimension $m$ if it is homeomorphic to Int $B^m$.
 The set $\dot{e}_i^m := \bar{e}_i^m - e_i^m$ for $m>0$ is called the \textbf{boundary} of the $m$-cell $e_i^m$.
\end{definition}

\begin{definition}\label{def:CWcomplex}
 A \textbf{CW complex} is a space $X$ and a collection of disjoint open cells $e_i$ whose union is $X$ such that:
 \begin{enumerate}
  \item $X$ is Hausdorff.
  \item For each open $m$-cell $e_i$ of the collection, there exist a continuous map $f_i: B^m \rightarrow X$ that maps Int$B^m$ homeomorphically onto $e_i$ and carries Bd$B^m$ into a finite union of open cells, each of dimension less than $m$.
  \item A set $A$ is closed in $X$ if $A\cap \bar{e}_i$ is closed in $\bar{e}_i$ for each $i$.
 \end{enumerate}
 If the collection of open cells is finite and the maps $f_i$ can be taken to be homeomorphisms, and each set $\dot{e}_i := \bar{e}_i - e_i$ equals the union of some open cells of $X$, then $X$ is called a \textbf{finite regular CW complex}.
\end{definition}

The condition (2) from the above definition is called ``closure-finiteness'' and condition (3) expresses the ``weak topology'' relative to the collection $\{\bar{e}_i\}$. These terms originate the letters C and W in the expression ``CW complex''.

In the following paragraphs we define a pseudomanifold as in \citep[Chapter IX]{Massey1991}. It is worth mentioning, that any regular CW complex on a closed connected $n$-manifold is an $n$-dimensional CW pseudomanifold.

\begin{definition}
We say that $e^m$ is a \textbf{face} of $e^n$ if $e^m \subset\bar{e}^n$, and denote as $e^m \preceq e^n$. If $e^m \neq e^n$, then $e^m$ is a \textbf{proper face} of $e^n$ ($e^m\prec e^n$).
The subspace $X^p$ of $X$ that is the union of the open cells of $X$ of dimension at most $p$ is called the \textbf{\textit{p}-skeleton} of $X$ and it is a CW complex in its own right.
\end{definition}

\begin{definition}\label{def:pseudomanifold}
 A \textbf{CW \textit{n}-pseudomanifold} is an $n$-dimensional finite regular CW complex which satisfies the following three conditions:
 \begin{enumerate}
  \item Every cell is a face of some $n$-cell.
  \item Every $(n-1)$-cell is a face of at most two $n$-cells.
  \item Given any two $n$-cells, $e^n_a$ and $e^n_b$, there exist a sequence of $n$-cells
  \[ e_0^n, e_1^n,\dots,e^n_k \]
  such that $e^n_a = e_0^n,\,e^n_b=e^n_k$, and $e^n_{i-1}$ and $e^n_i$ have a common $(n-1)$-dimensional face ($i=1,\dots,k$).
 \end{enumerate}
\end{definition}

If in the condition (2) of the above definition every $(n-1)$-dimensional cell is a face of exactly two $n$-cells, then we talk about \textbf{pseudomanifold without boundary}. On the other hand, if there exist at least one $(n-1)$-cell which is incident with only one $n$-cell, we talk about \textbf{pseudomanifold with boundary}.

The next definition is actually a result of a set of theorems in
\citep[Chapter IX]{Massey1991}, but for our purposes their implications are sufficient.

\begin{definition}
 Let $K=\{K^n\}$ be a pseudomanifold on the topological space $X$. And let $[ e_\lambda^n:e_\mu^{n-1}]$ be the \textbf{incidence number} of the cells $e_\lambda^n$ and $e_\mu^{n-1}$, for $n>0$, such that:
 \begin{enumerate}
  \item If $e_\mu^{n-1}$ is not a face of $e_\lambda^{n}$, then $[ e_\lambda^n:e_\mu^{n-1}]=0$.
  \item If $e_\mu^{n-1}$ is a face of $e_\lambda^{n}$, then $[ e_\lambda^n:e_\mu^{n-1}]=\pm 1$.
  \item If $e_\alpha^0$ and $e_\beta^{0}$ are the two vertices of the 1-cell $e_\lambda^{1}$, then $[e_\lambda^1:e_\alpha^{0}] + [e_\lambda^1:e_\beta^{0}]=0$.
  \item Let $e_\lambda^n$ and $e_\rho^{n-2}$ be cells such that $ e_\rho^{n-2}\prec e_\lambda^{n}$; let $e_\alpha^{n-1}$ and $e_\beta^{n-1}$ denote the unique $(n-1)$-cells $e^{n-1}$ such that $e_\rho^{n-2}\prec e^{n-1}\prec e_\lambda^{n}$. Then
  \[ [e_\lambda^n:e_\alpha^{n-1}][e_\alpha^{n-1}:e_\rho^{n-2}] + [e_\lambda^n:e_\beta^{n-1}][e_\beta^{n-1}:e_\rho^{n-2}] = 0. \]
 \end{enumerate}
With these conditions it is possible to choose the incidence numbers for each cell $e_\lambda^n$ in one and only one way. For $e^{n-1}\prec e^{n}$, we call the incidence number $[e^{n}:e^{n-1}]$ the \textbf{orientation} of $e^{n}$ with respect to $e^{n-1}$.
\end{definition}

Thus we can specify orientations for the cells by specifying a set of incidence numbers for the complex. We can now introduce the notion of orientability of a pseudomanifold.

\begin{definition}\label{def:orientablePseudomanifold}
 Let $K$ be an $n$-dimensional pseudomanifold (possibly with or without boundary), and let $e_1^n, e_2^n$ be $n$-cells with a common $(n-1)$-cell $e^{n-1}$.  We define orientations for $e_1^n$ and $e_2^n$ to be \textbf{coherent} (with respect to $e^{n-1}$) if:
 \[ [e_1^n: e^{n-1}]+[e_2^n: e^{n-1}] = 0. \]
 A set of orientations for all the $n$-cells of $K$ is said to be \textbf{coherent} if it is coherent in the above sense for any pair of $n$-cells with a common $(n-1)$-face.

 $K$ is said to be \textbf{orientable} if all its $n$-cells can be oriented such that any pair of $n$-cells sharing an $(n-1)$-dimensional face are oriented coherently. Otherwise it is called \textbf{nonorientable}.
\end{definition}

The orientability or nonorientability of an $n$-dimensional pseudomanifold $K$ only depends on the underlying topological space involved, and not on the choice of the regular CW complex $K$. We shall now define a chain on $K$, which together with the boundary operator give rise to the chain complex.

\begin{definition}\label{def:Chain}
 Let $K$ be a pseudomanifold. A \textbf{\textit{p}-chain} on $K$ is a function $c$ from the set of oriented $p$-cells of $K$ to the integers, such that $c(e) = -c({e}^\prime)$ if $e$ and ${e^\prime}$ are opposite orientations of the same cell. We add $p$-chains by adding their values, the resulting group is denoted $C_p(K)$ and is called the \textbf{chain group} of $K$. If $p<0$ or $p>\dim K$, we let $C_p(K)$ denote the trivial group.
\end{definition}

The incidence numbers, or equivalently, the set of orientations of $p$-cells of $K$, give us such a $p$-chain.

\begin{definition}\label{def:BoundaryOperator}
Let $C_n(K), n \geq 0,$ be the chain groups of $K$. The \textbf{boundary homomorphism}
\[ \partial_n: C_n(K) \rightarrow C_{n-1}(K),\; n>0, \]
is defined to be
\begin{equation}\label{eq:boundary}
 \partial_n(e^n) = \sum_\lambda [e^n:e^{n-1}_\lambda]e^{n-1}_\lambda,
\end{equation}
where $e^n$ is an \textit{oriented} $n$-cell of $K$ with $n>0$.
\end{definition}

From the definition of incidence numbers it follows that $\partial_n$ is well-defined and that $\partial_n(-e^n) = -\partial_n(e^n)$. The boundary homomorphism on 2--chains is illustrated in Figure \ref{fig1}.

\begin{definition}\label{def:ChainComplex}
 A \textbf{chain complex} $C_\ast=\{C_i,\partial_i\}$ is a sequence of abelian groups $C_i, i\in \mathbb{Z}$, and a sequence of homomorphisms $\partial_i:C_i\rightarrow C_{i-1}$ which are required to satisfy the condition
 \[ \partial_{i-1}\partial_i = 0 \;\; \forall i. \]
 For any such chain complex $C_\ast=\{C_i,\partial_i\}$ we have that $\im\partial_{i+1}\subset\ker\partial_i\subset C_i$ and we can define
 \[ H_i(C_\ast) = \frac{\ker\partial_i}{\im\partial_{i+1}}, \]
 called the $i$th \textbf{homology group of} $K$.
\end{definition}

If $K$ is an $n$-dimensional pseudomanifold, then $C_p(K) = 0$ for $p<0$ and $p>n$ and the chain complex reads:
\[0 \longrightarrow C_n(K) \overset{\partial_n}{\longrightarrow} \cdots \overset{\partial_{k+1}}{\longrightarrow} C_k(K) \overset{\partial_k}{\longrightarrow} \cdots \overset{\partial_1}{\longrightarrow} C_0(K) \longrightarrow 0.\]

\section{Discrete Differential Forms and the Exterior Derivative}\label{sec:DiscreteForms}
Having learned about chains and the boundary homomorphism we can now define cochains and the coboundary operator on them, which correspond to the \textit{discrete differential forms} and the \textit{discrete exterior derivative}, respectively.
\begin{definition}\label{def:cochain}
 Let $K$ be a pseudomanifold and $C_n(K)$ the group of oriented $n$-chains of $K$. Let $G$ be an abelian group. The group of \textbf{\textit{n}-dimensional cochains} of $K$, with coefficients in $G$, is the group
 \[C^n(K) =\text{ Hom }(C_n(K),G). \]
 The \textbf{coboundary operator} $d$ is defined to be the dual of the boundary operator $\partial_n: C_{n+1}(K) \rightarrow C_{n}(K)$, i.e., it is the homomorphism
 \[ d: C^n(K) \rightarrow C^{n+1}(K),\]
 such that $dd = 0.$
\end{definition}
We can think of a cochain group $C^n(K)$ as dual of a chain group $C_n(K)$. We now give the formal definition of discrete differential forms and the discrete exterior derivative:

\begin{definition}\label{def:DForm}
 A \textbf{discrete $G$--valued \textit{q}--form} $\alpha^q$ on a pseudomanifold $K$ is an element of $C^q(K)$, the group of \textit{q}--dimensional cochains of $K$, that is
 \[ \alpha^q \in C^q(K) =\text{ Hom }(C_q(K),G). \]
 The \textbf{discrete exterior derivative} $d_q: C^q(K) \rightarrow C^{q+1}(K)$ is the coboundary operator and it holds:
 \begin{equation}\label{eq:derivative}
  (d_q\alpha^q)(c_{q+1}) = \alpha(\partial_{q+1} c_{q+1}) = \sum_{c\in C_q(K)} [c_{q+1}:c]\alpha(c).
 \end{equation}
\end{definition}

\begin{figure}[ht]
\centering
\includegraphics[width=0.6\textwidth]{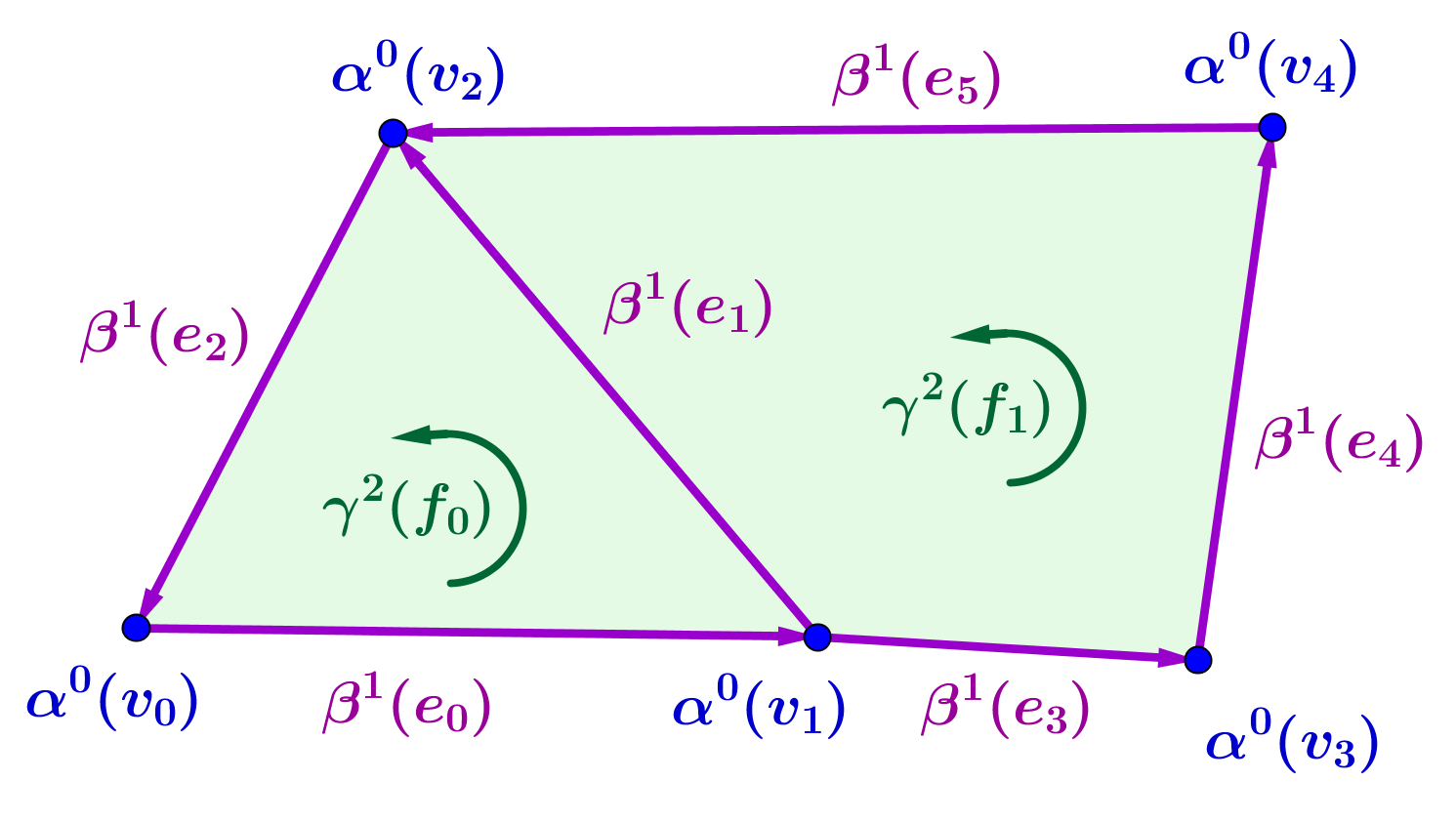}
\caption{
The boundary homomorphism and the exterior derivative on a polygonal complex. The pseudomanifold $K$ has 2--chain
$C_2 = (f_0, f_1)$, 1--chain $C_1 = (e_0,\dots,e_5$), and 0--chain $C_0 = (v_0,\dots,v_4)$. If 2--form $\gamma^2$ is given as the exterior derivative of 1--form $\beta^1$, then $\gamma = d\beta$ is defined per faces $f_0, f_1$ by equation (\ref{eq:derivative}). Concretely $\gamma(f_1) = d\beta(f_1) = \beta(\partial f_1)$ is:\\
$\beta(\partial f_1)
= \sum_{e_i \in C_1} [f_1:e_i]\beta(e_i)
= 0\cdot \beta(e_0) - 1\cdot \beta(e_1) + 0\cdot \beta(e_2)
+ 1\cdot \beta(e_3) + 1\cdot \beta(e_4) + 1\cdot \beta(e_5).
$
}\label{fig1}
\end{figure}

The discrete exterior derivative is illustrated in Figure \ref{fig1}. Since now on, we shall work with the coefficient group $G = (\mathbb{R},+)$, thus $C^\ast(K)$ are \textbf{real--valued discrete differential forms}.

Similarly to the definition of chain complexes (Definition \ref{def:ChainComplex}), we have the following:
\begin{definition}\label{def:CochainComplex}
A \textbf{cochain complex} $(C^\ast(K),d)$ on an $n$--dimensional pseudomanifold $K$ is a sequence of abelian groups of $q$-dimensional cochains $C^q(K)$ together with the coboundary operators
$d_q:C^q(K) \rightarrow C^{q+1}(K)$, i.e.,
\[0 \overset{d_n}{\longleftarrow} C^n(K) \overset{d_{n-1}}{\longleftarrow} \cdots \overset{d_{k}}{\longleftarrow} C^k(K) \overset{d_{k-1}}{\longleftarrow} \cdots \overset{d_0}{\longleftarrow} C^0(K) \longleftarrow 0.\]

 A $q$--form $\alpha$ is said to be a \textbf{closed form} (cocycle) if $d\alpha = 0$, and it is called an \textbf{exact form} (coboundary) if there exist a $(q-1)$--form $\beta$ such that $d\beta =\alpha$.
 
Note that each exact $q$-form $\alpha = d\beta$ is closed, since $d\alpha = dd\beta=0$. Two closed $q$-forms are \textbf{cohomologous} if they differ by an exact $q$--form. We denote cohomology classes by $H^q$:
 \[H^q(C^\ast)=\frac{\ker d_q}{\im d_{q-1}}.\]
\end{definition}

\medskip
Thus a discrete $q$--form $\alpha^q$ is a $q$--cochain defined per each element of a $q$--chain $C_q(K)$:
\[
 \alpha^q = \alpha(c_0),\dots,\alpha(c_p) ,\;\; c_i\in C_q(K).
\]
We sum two discrete $q$--forms $\alpha^q,\beta^q$ by summing their values per each element of a $q$--chain $C_q(K)$ and multiply a $q$--form $\alpha^q$ by a scalar $t\in\Re$ in element--wise fashion too:
\begin{align}
\alpha^q + \beta^q &= (\alpha +\beta)(c_0),\dots,(\alpha +\beta)(c_p), \label{eq:polygonalSum}\\
t\cdot \alpha^q &= t\cdot \alpha(c_0),\dots,t\cdot \alpha(c_p).\label{eq:polygonalScalarMultiplication}
\end{align}


\section{Wedge and Cup Products}\label{sec:CupProduct}
In this section we first present the wedge product on a manifold and then give a general definition of an abstract cup product of cellular cochains to show the clear analogy between these two. We then proceed with equations for cup products on triangles and quadrilaterals, and finally present simple formulas for cup products on general polygons in Section \ref{subsec:CupOnPolygons}.

\subsubsection*{Wedge product}
A differential $k$-form $\alpha$ on a smooth $n$--manifold $M$ is a tensor field of type $(0,k)$ that is completely antisymmetric, i.e., for $\alpha \in \Omega^k(M) \subset T^0_k(M)$ we have $\alpha(v,w) = -\alpha(w,v)$ for all $v,w\in TM$. Let $\{e_1,\dots,e_n\}$ be a basis of the tangent space $T_x M$ to $M$ at point $x\in M$. Let $\alpha \in T^0_k(M)$ and $\beta \in T^0_l(M)$, we define their \textbf{wedge product} $\alpha\wedge\beta\in T^0_{k+l}(M)$ at the point $x\in M$ by
\begin{equation}\label{eq:WedgeAbraham1988}
 (\alpha\wedge\beta) = \sum_{{\tau}}\sign{\tau} \alpha(e_{\tau(1)},\dots,e_{\tau(k)})\beta(e_{\tau(k+1)},\dots,e_{\tau(k+l)}),
\end{equation}
where the sum is over all permutations ${\tau}$ of $\{1,\dots,k+l\}$ such that $\tau(1)<\dots<\tau(k)$ and $\tau(k+1)<\dots<\tau(k+l)$.
And we have that $\sign{\tau} = +1$ if the permutation is even and $\sign{\tau} = -1$ if it is odd.

The wedge product exhibits the following properties, which we maintain also in the discrete setting except the associativity (compare with the later Proposition \ref{prop:PolygonalCup}):

\begin{proposition}\label{prop:wedgeProduct}
For $\alpha \in T^0_k(M)$, $\beta \in T^0_l(M)$, and $\gamma \in T^0_m(M)$, we have
\begin{enumerate}
 \item $\wedge$ is bilinear: $\alpha\wedge(c_1\beta+c_2\gamma) = c_1(\alpha\wedge\beta)+c_2(\alpha\wedge\gamma)$
 and $(c_1\alpha+c_2\beta)\wedge\gamma = c_1(\alpha\wedge\gamma)+c_2(\beta\wedge\gamma)$ for some constants $c_1,c_2\in\mathbb{R}$,
 \item $\wedge$ is skew commutative: $\alpha\wedge\beta = (-1)^{kl}\beta\wedge\alpha$,
 \item $\wedge$ is associative: $\alpha\wedge(\beta\wedge\gamma) = (\alpha\wedge\beta)\wedge\gamma$,
 \item Leibniz rule: $d (\alpha\wedge\beta) = d\alpha\wedge\beta + (-1)^k\alpha\wedge d\beta$.
\end{enumerate}
\end{proposition}
For the proof see \citep[Proposition 7.1.5 and Theorem 7.4.1]{Abraham1988}.

\subsubsection*{Cup product}
The wedge product allows for building higher degree forms from lower degree ones, similarly a cup product is a product of cochains of arbitrary degree $p$ and $q$ that returns a cochain of degree $p+q$.

Whitney in \citep[\S 9 of Appendix II]{Whitney1957} gives the following abstract definition of a cup product, which also appears in \citep[Definition 2.3.3]{Arnold2012}.

\begin{definition}\label{def:GeneralCupProduct}
Let $X$ be a cell complex, the \textbf{cup product} of two cochains $c^p$ and $c^q$ is a bilinear operation
$\cup : C^p(X)\times C^q(X)\rightarrow C^{p+q}(X)$ that satisfies the following three properties:
\begin{enumerate}
 \item Let $\sigma_p\in C_p(X)$ and $\sigma_q\in C_q(X)$. Then $\sigma^p\cup\sigma^q$ is a $(p+q)$-cochain in $St(\sigma_p)\cdot St(\sigma_q)$,
 where $St(\sigma_i)$ is the union of all cells in which $c_i$ is a face, and $A\cdot B$ denotes the union of all cells in $A$ and $B$.
 \item $d (c^p\cup c^q) = d c^p\cup c^q + (-1)^pc^p\cup d c^q$ (Leibniz rule).
 \item If $X$ is connected, then there exist a real number $\gamma_\cup$ such that
 $I^0\cup c^p = c^p \cup I^0 = \gamma_\cup c^p$, where $I^0$ is the constant 0-cochain that takes value 1 on the 0-cells of $X$.
\end{enumerate}
\end{definition}

The interested reader is also referred to \citep[Chapter XIII]{Massey1991}. Whitney further asserts that thanks to the Leibniz rule being valid, we have the following:
\begin{proposition}\label{prop:GeneralCup}
For $\alpha \in H^k(X)$, $\beta \in H^l(X)$, $\gamma\in H^m(X)$, the cup product defines a bilinear operation $H^k(X)\times H^l(X)\rightarrow H^{k+l}(X)$, which is uniquely determined (well--defined) and is
\begin{enumerate}
 \item skew commutative: $\alpha\cup\beta = (-1)^{kl}\beta\cup\alpha$,
 \item associative: $\alpha\cup(\beta\cup\gamma) = (\alpha\cup\beta)\cup\gamma$.
\end{enumerate}
Thus the cup endows the cellular cohomology $H^*(X)$ with a \textbf{graded commutative ring} structure, with the cohomology class of $I^0$ as unit element.
\end{proposition}

Comparing Proposition \ref{prop:wedgeProduct} with Definition \ref{def:GeneralCupProduct} and Proposition \ref{prop:GeneralCup}, we affirm that both the wedge and the cup product, respectively, are bilinear operations that take as input two forms, resp. cochains, of arbitrary degree $p$ and $q$ and return a form, resp. cochain, of degree $p+q$. Moreover, they both satisfy the Leibniz rule. But a general cup products is guaranteed to be associative and skew--commutative only on cohomology. Fortunately, for our cup product on polygons we will recuperate the skew--commutativity for any discrete forms (not just closed ones). However, our polygonal cup product will be associative only on closed forms, see the Proposition \ref{prop:PolygonalCup}.

Even though explicit formulas for computing cup product on a general cell complex are unknown, there are well--known explicit formulas for cup products on simplicial and cubical complexes \citep{Arnold2012}. We expand the set of known cup products with \textit{cup products on general polygons with coefficients in $\mathbb{R}$, which is one of our contributions}.

\subsection{Simplicial and cubical cup products}
Let us first recall the definition of \textbf{simplicial Whitney forms}. Let $c = (v_0,\dots,v_p)$ be a $p$--simplex with barycentric coordinates $x_0,\dots,x_p$. Thus $0\leq x_j\leq 1\; \forall j$ and $\sum_{i=0}^p = 1$. The simplicial Whitney $p$--form corresponding to the simplex $c$ is
\[
W(c) = p! \sum_{i=0}^p (-1)^i x_i dx_0\wedge\dots \wedge \widehat{dx_i} \wedge\dots\wedge dx_p,
\]
where $\;{}\widehat{}\;$ indicates a term omitted from the product.
This definition extends linearly to define Whitney forms on all of $C^\ast(K)$.

Arnold in \citep[Definition 4.2.1]{Arnold2012} presents the following definition of a cup product on simplices, which she calls the \textbf{Wilson's cup product}, $\cup: C^p(K)\times C^q(K)\rightarrow C^{p+q}(K)$ given by
\begin{equation}\label{eq:WilsonsCup}
 \alpha\cup\beta (c) = \int_c W\alpha\wedge W\beta,
\end{equation}
where $c$ is an $(p+q)$--simplex and $W$ is the simplicial Whitney form.

The above cup product was given already in \citep{Whitney1957}. Scott Wilson later proved that the cup product converges to wedge product of differential forms \citep[Theorem 5.4]{Wilson2007}. Furthermore, Arnold claims that Wilson's work was the motivation for the majority of results presented in \citep{Arnold2012}, including the definition of the cubical cup product.

\begin{figure}[h]
\centering\includegraphics[width=0.95\textwidth]{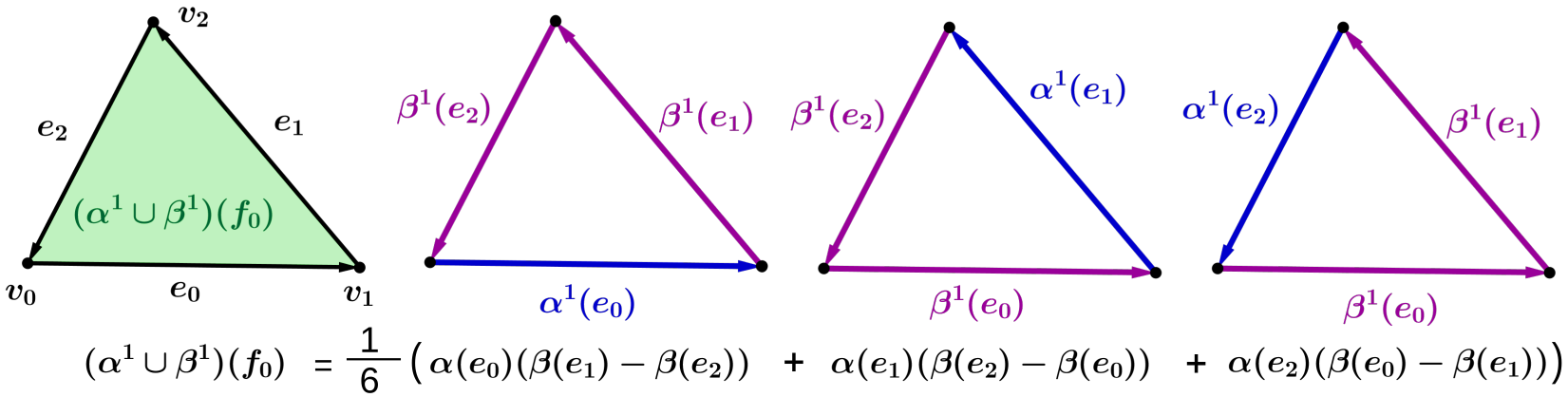}
\caption{The cup product of two 1--forms is a 2--form located on faces.}\label{fig2}
\end{figure}

Using the equation (\ref{eq:WilsonsCup}) we compute cup products on 2--dimensional simplicial complexes and state them as the next definition. The cup product on 1--cochains is illustrated in Figure \ref{fig2}.

\begin{definition}\label{def:SimplicialCup}
Let $K$ be a 2--pseudomanifold whose 2--cells are triangles. The \textbf{simplicial cup product} $\cup: C^p(K,\Re)\times C^q(K,\Re)\rightarrow C^{p+q}(K,\Re)$ is defined by its action on a $(p+q)$--dimensional simplex $(v_0,\dots,v_{(p+q)})$ as:
\begin{eqnarray*}
(\alpha^0\cup\beta^0)(v_0)  &=&  \alpha(v_0)\beta(v_0),\\
(\alpha^0\cup\beta^1)(v_0,v_1)  &=&  \frac{1}{2}(\alpha(v_0)+\alpha(v_1))\beta(v_0,v_1),\\
(\alpha^0\cup\beta^2)(v_0,v_1,v_2) &=&\frac{1}{3}(\alpha(v_0)+\alpha(v_1)+\alpha(v_2))\beta(v_0,v_1,v_2),  \\
(\alpha^1\cup\beta^1)(v_0,v_1,v_2) &=& \frac{1}{6}\sum_{i=0}^2 \alpha(v_i,v_{i+1})(\beta(v_{i+1},v_{i+2}) - \beta(v_{i-1},v_{i})),
\end{eqnarray*}
where the indices are taken to be modulo 3.
\end{definition}

\medskip
In \citep[Definition 3.2.1]{Arnold2012} the author defines a cubical cup product $\cup_c$, for which the product of two cochains agrees with the wedge product of their \textbf{cubical Whitney forms} \citep[Section 3.2.2]{Arnold2012}, thus the Whitney map provides a connection of cubical cohomology with de Rham cohomology \citep[Theorem 3.2.12]{Arnold2012}. From \citep[Definition 3.2.1]{Arnold2012} we computed the cup product on $k$--dimensional cubes for $0\leq k\leq 2$ and we present the formulas in the following definition:

\begin{definition}\label{def:CubicalCup}
Let $K$ be a 2--pseudomanifold whose 2--cells are quadrilaterals.
The \textbf{cubical cup product} $\cup: C^p(K,\Re)\times C^q(K,\Re)\rightarrow C^{p+q}(K,\Re)$ is defined by its action on a $(p+q)$--dimensional cube $(v_0,\dots,v_{2^{(p+q)}-1})$ as:
\begin{eqnarray*}
(\alpha^0\cup\beta^0)(v_0)  &=&  \alpha(v_0)\beta(v_0),\\
(\alpha^0\cup\beta^1)(v_0,v_1)  &=&  \frac{1}{2}(\alpha(v_0)+\alpha(v_1))\beta(v_0,v_1),\\
(\alpha^0\cup\beta^2)(v_0,v_1,v_2,v_3) &=&\frac{1}{4}\Big(\sum_{i=0}^3\alpha(v_i)\Big)\beta(v_0,v_1,v_2,v_3),  \\
(\alpha^1\cup\beta^1)(v_0,v_1,v_2,v_3) &=& \frac{1}{4}\sum_{i=0}^3 \alpha(v_i,v_{i+1})(\beta(v_{i+1},v_{i+2}) - \beta(v_{i-1},v_{i})),
\end{eqnarray*}
where the indices are taken to be modulo 4.
\end{definition}

\subsection{The cup product on polygons}\label{subsec:CupOnPolygons}
As advertised earlier, we can extend the definitions of the discrete wedge product to any $p$--polygon such that it still satisfies the defining properties of a cup product (Definition \ref{def:GeneralCupProduct}) and therefore automatically also the properties of Proposition \ref{prop:GeneralCup}. Moreover, we show that our polygonal cup product is a skew--commutative bilinear operation on any polygonal forms, not just on closed ones, see the Proposition \ref{prop:PolygonalCup}.

\begin{definition}\label{def:PolygonalCup}
Let $K$ be a 2--pseudomanifold whose 2--cells are simple polygons. The cup product $C^k(K,\Re)\times C^l(K,\Re)\rightarrow C^{k+l}(K,\Re)$ of two discrete forms $\alpha^k,\beta^l$ is a $(k+l)$--form $\alpha^k\cup\beta^l$ defined on each $(k+l)$--cell of $C_{k+l}(K)$. Let $f=(v_0,\dots,v_{p-1})$ be a $p$--polygonal face,
$e_i = (v_i,v_{i+1})$ an edge, and $v$ a vertex of $K$, the \textbf{polygonal cup product} is defined for each degree and per each $(k+l)$--cell as:
\begin{align}
(\alpha^0\cup\beta^0)(v)  &=  \alpha(v)\beta(v), \label{eq:PolygonalCup00} \\
(\alpha^0\cup\beta^1)(e_i)  &=  \frac{1}{2}(\alpha(v_i)+\alpha(v_{i+1}))\beta(e_i),\label{eq:PolygonalCup01} \\
(\alpha^0\cup\beta^2)(f) &=\frac{1}{p}\bigg( \sum_{i=0}^{p-1}\alpha(v_i) \bigg) \beta(f),\label{eq:PolygonalCup02}   \\
(\alpha^1\cup\beta^1)(f) &= \sum_{i=0}^{p-1} \alpha(i)
\sum_{a=1}^{\lfloor\frac{p-1}{2} \rfloor} \bigg( \frac{1}{2} - \frac{a}{p} \bigg) (\beta((i+a)\%p) - \beta((i-a)\%p)),\label{eq:PolygonalCup11}
\end{align}
where $\alpha(i):=\alpha(e_i)=\alpha(v_i,v_{(i+1)\%p})$ and $\%p$ means modulo $p$.
\end{definition}

It is easy to see that the simplicial and cubical cup products (Definitions \ref{def:SimplicialCup} and \ref{def:CubicalCup}) are just special cases of the product defined above.

In Theorem \ref{theo:PolygonalCupIsCup} we prove that the polygonal cup product actually \textit{is} a cup product, that is, we prove that it satisfies the defining properties of a cup product (Definition \ref{def:GeneralCupProduct}). But we first prove, in the following lemma, that it obeys the Leibniz product rule with the discrete derivative.

\begin{lemma}\label{lemma:LeibnizPolygonalCup}
 The discrete wedge product from Definition \ref{def:PolygonalCup} satisfies the Leibniz rule with the discrete exterior derivative, that is:
 \begin{align*}
d_0(\alpha^0\cup\beta^0) &=  d_0\alpha\cup\beta + \alpha\cup d_0\beta,\\
d_1(\alpha^0\cup\beta^1) &=  d_0\alpha\cup\beta + \alpha\cup d_1\beta.
\end{align*}
For $d(\alpha^0\cup\beta^2)$ and $d(\alpha^1\cup\beta^1)$ we get trivially 0.
\end{lemma}

\begin{proof} We will prove the Leibniz rule for each degree separately. For the cases $d(\alpha^0\cup\beta^2)$ and $d(\alpha^1\cup\beta^1)$ we get trivially 0 because there are no 3--dimensional chains (and thus no 3--forms) on a 2--dimensional pseudomanifold.\\
\textbf{The case $d(\alpha^0\cup\beta^0)$:} let $e=(v_i,v_j)$ be an edge of $K$, then
\begin{align*}
 (d_0\alpha\cup\beta + \alpha\cup d_0\beta)(e) &= (\alpha(v_{j})-\alpha(v_i))\frac{\beta(v_i)+\beta(v_{j})}{2}\\
 &+ \frac{\alpha(v_i)+\alpha(v_{j})}{2}(\beta(v_{j})-\beta(v_i))\\
 & = \frac{\beta(v_i)(-2\alpha(v_i)) + \beta(v_{j})(2\alpha(v_{j}))}{2}\\
 & = \alpha(v_{j})\beta(v_{j}) - \alpha(v_i)\beta(v_i) \\
 & = (\alpha\cup\beta)(v_{j}) - (\alpha\cup\beta)(v_i)\\
 & = d_0(\alpha^0\cup\beta^0)(e).
\end{align*}
\\
\textbf{The case $d(\alpha^0\cup\beta^1)$:} Let us first simplify the notation and denote $\alpha_i = \alpha^0(v_i)$, $\beta_i=\beta^1(e_i)$. Let also any index $i$ to be understood as $i\%n$. Let $f=(v_0,\dots,v_{n-1})$ be an $n$--polygonal face with edges $e_i=(v_i,v_{i+1})$. Due to the indexes being cyclic (modulo $n$), we can reorder the following sums in this manner:
\begin{align*}
 d(\alpha^0\cup\beta^1)(f) & = \sum_{i=0}^{n-1} \frac{\alpha_i + \alpha_{i+1}}{2}\beta_i = \sum_{i=0}^{n-1} \frac{\beta_{i-1} + \beta_{i}}{2}\alpha_i\\
 (d\alpha^0\cup\beta^1)(f) & = \sum_{i=0}^{n-1} \sum_{a=1}^{\lfloor\frac{n-1}{2} \rfloor} \bigg(\frac{1}{2}-\frac{a}{n}\bigg)(\alpha_{i+1} - \alpha_i)(\beta_{i+a} - \beta_{i-a}) \\
 & = \sum_{i=0}^{n-1} \sum_{a=1}^{\lfloor\frac{n-1}{2} \rfloor} \frac{n-2a}{2n}(\alpha_i)(\beta_{i+a-1} - \beta_{i-a-1} - \beta_{i+a} + \beta_{i-a})\\
 (\alpha^0\cup d\beta^1)(f)& = \frac{1}{n}\sum_{i=0}^{n-1}\alpha_i \sum_{i=0}^{n-1}\beta_i
\end{align*}
Thus we want to show that
\[
 \sum_{i=0}^{n-1} \alpha_i \frac{\beta_{i-1} + \beta_{i}}{2}
 = \sum_{i=0}^{n-1}\frac{\alpha_i}{2n}\bigg[ \sum_{a=1}^{\lfloor\frac{n-1}{2} \rfloor} (n-2a)(\beta_{i+a-1} - \beta_{i+a} + \beta_{i-a} - \beta_{i-a-1})
+ 2\sum_{j=0}^{n-1}\beta_j\bigg],
\]
but instead we will show that the following equation (\ref{eq:Alpha0ToExpandLeibniz}) holds for any $i\in\{0,\dots,n-1\}$, which implies that also the equality above is true.
\begin{equation}\label{eq:Alpha0ToExpandLeibniz}
\frac{\beta_{i-1} + \beta_{i}}{2} = \frac{1}{2n}\bigg[\sum_{a=1}^{\lfloor\frac{n-1}{2} \rfloor} (n-2a)(\beta_{i+a-1} - \beta_{i+a} + \beta_{i-a} - \beta_{i-a-1}) + 2\sum_{j=0}^{n-1}\beta_j\bigg]
\end{equation}

If $n$ is even, then $a=1,\dots,\lfloor\frac{n-1}{2} \rfloor$ is equivalent to $a = 1,2,\dots,\frac{n}{2}-1$, and if $n$ is odd, then it is equivalent to $a= 1,\dots,\frac{n-1}{2}$. We will first show the validity of (\ref{eq:Alpha0ToExpandLeibniz}) for $n$ even and then for $n$ odd.

\textbf{For $n$ even} we compute:
\begin{align*}
\sum_{a=1}^{\frac{n}{2}-1}(n-2a)&(\beta_{k+a-1} - \beta_{k+a}) = (n-2)(\beta_{k} - \beta_{k+1})+(n-4)(\beta_{k+1} - \beta_{k+2})\\
  &\;+ (n-6)(\beta_{k+2} - \beta_{k+3}) + \cdots + 6(\beta_{k+\frac{n}{2}-4} - \beta_{k+\frac{n}{2}-3})\\
  &\;+ 4(\beta_{k+\frac{n}{2}-3} - \beta_{k+\frac{n}{2}-2}) + 2(\beta_{k+\frac{n}{2}-2} - \beta_{k+\frac{n}{2}-1})\\
  =&\; n\beta_k - 2(\beta_k + \beta_{k+1} + \beta_{k+2} +  \beta_{k+3} +\dots\\
  & + \beta_{k+\frac{n}{2}-3} + \beta_{k+\frac{n}{2}-2} + \beta_{k+\frac{n}{2}-1}),\\
\sum_{a=1}^{\frac{n}{2}-1}(n-2a)&(\beta_{k-a} - \beta_{k-a-1}) = (n-2)(\beta_{k-1} - \beta_{k-2})+(n-4)(\beta_{k-2} - \beta_{k-3})\\
  &\;+ (n-6)(\beta_{k-3} - \beta_{k-4}) + \cdots + 6(\beta_{k-\frac{n}{2}+3} - \beta_{k-\frac{n}{2}+2}) \\
  &\;+ 4(\beta_{k-\frac{n}{2}+2} - \beta_{k-\frac{n}{2}+1}) + 2(\beta_{k-\frac{n}{2}+1} - \beta_{k-\frac{n}{2}})\\
  =&\; n\beta_{k-1} - 2(\beta_{k-1} + \beta_{k-2} +  \beta_{k-3} +
  \dots\\
  &+ \beta_{k+\frac{n}{2}+2} + \beta_{k+\frac{n}{2}+1} + \beta_{k+\frac{n}{2}}).
\end{align*}
Thus
\begin{align*}
 \sum_{a=1}^{\frac{n}{2}-1}(n-2a)(\beta_{k+a-1}-\beta_{k+a}+\beta_{k-a} -\beta_{k-a-1} ) &= n(\beta_{k-1}+\beta_k) -2\sum_{i=k}^{k+n-1}\beta_i \\
 &= n(\beta_{k-1}+\beta_k) -2\sum_{j=0}^{n-1}\beta_j
\end{align*}
and for $n$ even we obtain
\[
\frac{1}{2n}\bigg[\sum_{a=1}^{\lfloor\frac{n-1}{2} \rfloor}(n-2a)(\beta_{k+a-1}-\beta_{k+a}+\beta_{k-a} -\beta_{k-a-1} ) + 2\sum_{j=0}^{n-1}\beta_j\bigg] = \frac{\beta_{k-1}+\beta_k}{2}.
\]

\textbf{For $n$ odd} we proceed in a similar fashion:
\begin{align*}
\sum_{a=1}^{\frac{n-1}{2}}(n-2a)&(\beta_{k+a-1} - \beta_{k+a}) = (n-2)(\beta_{k} - \beta_{k+1})+(n-4)(\beta_{k+1} - \beta_{k+2})\\
  &\;+ (n-6)(\beta_{k+2} - \beta_{k+3}) + \cdots + 5(\beta_{k+\frac{n-1}{2}-3} - \beta_{k+\frac{n-1}{2}-2})\\
  &\;+ 3(\beta_{k+\frac{n-1}{2}-2} - \beta_{k+\frac{n-1}{2}-1}) + 2(\beta_{k+\frac{n-1}{2}-1} - \beta_{k+\frac{n-1}{2}})\\
  =&\; n\beta_k - 2(\beta_k + \beta_{k+1} + \beta_{k+2} +  \beta_{k+3} +\dots+ \beta_{k+\frac{n-1}{2}-2} + \beta_{k+\frac{n-1}{2}-1}) \\
  &- \beta_{k+\frac{n-1}{2}},\\
\sum_{a=1}^{\frac{n-1}{2}}(n-2a)&(\beta_{k-a} - \beta_{k-a-1}) = (n-2)(\beta_{k-1} - \beta_{k-2})+(n-4)(\beta_{k-2} - \beta_{k-3})\\
  &\;+ (n-6)(\beta_{k-3} - \beta_{k-4}) + \cdots + 5(\beta_{k-\frac{n-1}{2}+2} - \beta_{k-\frac{n-1}{2}+1}) \\
  &\;+ 3(\beta_{k-\frac{n-1}{2}+1} - \beta_{k-\frac{n-1}{2}}) + (\beta_{k-\frac{n-1}{2}} - \beta_{k-\frac{n-1}{2}-1})\\
  =&\; n\beta_{k-1} - 2(\beta_{k+\frac{n-1}{2}+1} + \beta_{k+\frac{n-1}{2}+2} +\dots+ \beta_{k-3} + \beta_{k-2} +  \beta_{k-1})\\
  &- \beta_{k+\frac{n-1}{2}}.
\end{align*}
Thus
\begin{align*}
 \sum_{a=1}^{\frac{n-1}{2}}(n-2a)(\beta_{k+a-1}-\beta_{k+a}+\beta_{k-a} -\beta_{k-a-1} ) &= n(\beta_{k-1}+\beta_k) -2\sum_{i=k}^{k+n-1}\beta_i \\
 &= n(\beta_{k-1}+\beta_k) -2\sum_{j=0}^{n-1}\beta_j
\end{align*}
and for $n$ odd we get
\[
\frac{1}{2n}\bigg[\sum_{a=1}^{\lfloor\frac{n-1}{2} \rfloor}(n-2a)(\beta_{k+a-1}-\beta_{k+a}+\beta_{k-a} -\beta_{k-a-1} ) + 2\sum_{j=0}^{n-1}\beta_j\bigg] = \frac{\beta_{k-1}+\beta_k}{2}.
\]
We have shown that the equality (\ref{eq:Alpha0ToExpandLeibniz}) holds for any $i\in\{ 0,1,\dots,n-1\}$, which implies $d_1(\alpha^0\cup\beta^1)=d_0\alpha\cup\beta + \alpha\cup d_1\beta$.
\end{proof}
%
%

Now we are ready to claim that our discrete wedge product on polygonal pseudomanifolds is a cup product:

\begin{theorem}\label{theo:PolygonalCupIsCup}
The polygonal cup product from Definition \ref{def:PolygonalCup} is a cup product, i.e., it satisfies the Definition \ref{def:GeneralCupProduct}.
\end{theorem}

\begin{proof}
From the formulas (\ref{eq:PolygonalCup00})--(\ref{eq:PolygonalCup11}) and the way we sum two discrete forms (eq. \ref{eq:polygonalSum}) or multiply them by scalar (eq. \ref{eq:polygonalScalarMultiplication}), it is a simple exercise to see that
$\forall t\in\Re,\; \forall c_i\in C_{k+l}, 0\leq (k+l) \leq 2$:
\begin{align*}
(t\cdot\alpha^k\cup\beta^l)(c_i)  &= (\alpha^k\cup t\cdot\beta^l)(c_i),\\
((\alpha^k_1 + \alpha^k_2) \cup\beta^l)(c_i)  &= (\alpha^k_1 \cup\beta^l)(c_i) +
(\alpha^k_2 \cup\beta^l)(c_i),\\
(\alpha^k \cup (\beta^l_1 + \beta^l_2)(c_i)  &= (\alpha^k \cup\beta^l_1)(c_i) +
(\alpha^k \cup\beta^l_2)(c_i).
\end{align*}
Therefore the polygonal cup product is a \textit{bilinear} operation.

It also satisfies item 1 of Definition \ref{def:GeneralCupProduct} -- the result of product of two forms is located on a face that is incident to both the elements where the two original discrete forms reside.

Due to Lemma \ref{lemma:LeibnizPolygonalCup}, the product satisfies item 2, the Leibniz rule.

We will now prove that item 3 of Definition \ref{def:GeneralCupProduct} also holds. A pseudomanifold $K$ as we defined it is a connected cell complex, therefore we have to show that there exist a real number $\xi$ such that $I^0\cup \alpha^p = \alpha^p\cup I^0 = \xi\alpha^p$, where $I^0$ is the constant 0--form that takes value 1 on the vertices of $K$. Using the equations of Definition \ref{def:PolygonalCup}, we can thus write:
 \begin{align*}
  I^0\cup\alpha^0 &= \sum_{v_i\in V(K)}I^0\cup\alpha^0(v_i) =  \sum_{v_i\in V}I^0(v_i)\alpha^0(v_i) = \sum_{v_i\in V}1\alpha^0(v_i) = 1\alpha^0,\\
  \alpha^0\cup I^0 &= \sum_{v_i\in V}\alpha^0\cup I^0(v_i) = \sum_{v_i\in V}\alpha^0(v_i)I^0(v_i) = \sum_{v_i\in V}\alpha^0(v_i)1 = 1\alpha^0,\\
  I^0\cup\alpha^1 &= \sum_{(v_i,v_j)\in E(K)}I^0\cup\alpha^1(v_i,v_j) = \sum_{(v_i,v_j)\in E}\frac{I^0(v_i)+I^0(v_j)}{2}\alpha^1(v_i,v_j) \\
  &= \sum_{(v_i,v_j)\in E}\frac{1+1}{2}\alpha^1(v_i,v_j) = 1\alpha^1,\\
  \alpha^1\cup I^0 &= \sum_{(v_i,v_j)\in E}\alpha^1\cup I^0(v_i,v_j) = \sum_{(v_i,v_j)\in E}\alpha^1(v_i,v_j)\frac{I^0(v_i)+I^0(v_j)}{2} \\
  &= \sum_{(v_i,v_j)\in E}\alpha^1(v_i,v_j) \frac{1+1}{2}= 1\alpha^1,\\
  I^0\cup\alpha^2 &= \sum_{f_i\in F(K)}I^0\cup\alpha^2(f_i) = \sum_{f_i\in F}\sum_{v_j\prec f_i}\frac{I^0(v_j)}{p_i}\alpha^2(f_i) \\
  &= \sum_{f_i\in F}\frac{p_i}{p_i}\alpha^2(f_i) = 1\alpha^2,\\
  \alpha^2\cup I^0 &= \sum_{f_i\in F}\alpha^2\cup I^0(f_i) = \sum_{f_i\in F}\alpha^2(f_i)\sum_{v_j\prec f_i}\frac{I^0(v_j)}{p_i} \\
  &= \sum_{f_i\in F}\alpha^2(f_1)\frac{p_i}{p_i}= 1\alpha^2,
 \end{align*}
where $p_i$ is the number of vertices of a face $f_i$. We thus conclude that $\xi=1$ satisfies $I^0\cup \alpha^p = \alpha^p\cup I^0 = \xi\alpha^p \;\; \forall p\in\{0,1,2\}$.
\end{proof}

In the following proposition, which can be thought of as an analog of Proposition \ref{prop:wedgeProduct}, we summarize and prove further properties of our polygonal cup product.

\begin{proposition}\label{prop:PolygonalCup}
Let $K$ be a 2--pseudomanifold, $\alpha^k \in C^k(K)$, $\beta^l \in C^l(K)$, and $\gamma^m \in C^m(K)$.
The polygonal cup product from Definition \ref{def:PolygonalCup} satisfies these properties:
\begin{enumerate}
 \item Bilinearity.
 \item Skew--commutativity $\alpha^k \cup\beta^l = (-1)^{kl} \beta^l\cup\alpha^k$.
 \item Leibniz rule.
 \item Associativity under the following conditions:\\
    $\alpha^0\cup(\beta^0\cup\gamma^0)=(\alpha^0\cup\beta^0)\cup\gamma^0$ always,\\
    $\alpha^0\cup(\beta^0\cup\gamma^1)=(\alpha^0\cup\beta^0)\cup\gamma^1$ if $d\alpha^0 = 0$ or $d\beta^0 = 0$,\\
    $\alpha^0\cup(\beta^0\cup\gamma^2)=(\alpha^0\cup\beta^0)\cup\gamma^2$ if $d\alpha^0 = 0$ or $d\beta^0 = 0$,\\
    $\alpha^0\cup(\beta^1\cup\gamma^1)=(\alpha^0\cup\beta^1)\cup\gamma^1$ if $d\alpha^0 = 0$.
\item The error of associativity vanishes under mesh refinement: \\
$\big(\alpha^0\cup(\beta^k\cup\gamma^l) - (\alpha^0\cup\beta^k)\cup\gamma^l\big)(c_{k+l})\to 0$ for $\norm{v_i-v_j}\to 0$ $\forall v_i,v_j\in c_{k+l}$
if $d\alpha,d\beta,d\gamma$ are bounded.
\end{enumerate}
\end{proposition}

\begin{proof}
 From Theorem \ref{theo:PolygonalCupIsCup} we know the product is a \textit{bilinear} operation satisfying the \textit{Leibniz rule} with discrete exterior derivative. We now show that it is also skew--commutative on any discrete forms and associative in the declared fashion.

Concerning the \textit{skew--commutativity}, from the equations (\ref{eq:PolygonalCup00})--(\ref{eq:PolygonalCup02}) it is easy to see that $\alpha^0 \cup\beta^l = \beta^l\cup\alpha^0$ for $l=0,1,2$. The case of the product of two 1--forms is slightly more involving. We need to prove that $(\alpha^1\cup\beta^1)(f) = - (\beta^1\cup\alpha^1)(f)$ on an $n$--polygonal face $f=(v_0,\dots,v_{n-1})$. Let us remember that
\[
 (\alpha^1\cup\beta^1)(f) = \sum_{a=1}^{\lfloor\frac{n-1}{2} \rfloor} \bigg( \frac{1}{2} - \frac{a}{n} \bigg) \sum_{i=0}^{n-1} \alpha(i)(\beta((i+a)\%n) - \beta((i-a)\%n),
\]
where $\alpha(i):=\alpha(e_i)=\alpha(v_i,v_{(i+1)\%n})$ and $\%n$ means modulo $n$. Again, to simplify the notation we will omit writing the symbol $\%n$ in the indexes and by some index $k$ we will always understand $k\%n$. We next show that
\[
\sum_{i=0}^{n-1} \alpha(i)(\beta(i+a) - \beta(i-a)) = -\sum_{i=0}^{n-1} \beta(i)(\alpha(i+a) - \alpha(i-a)).
\]
Due to the indexes being cyclic (modulo $n$), we can write:
\small
\begin{align*}
 \sum_{i=0}^{n-1} \alpha(i)\beta(i+a) &= \alpha(0)\beta(a) + \alpha(1)\beta(1+a) + \cdots + \alpha(n-a-1)\beta(n-1)\\
 &+ \alpha(-a)\beta(0) + \alpha(1-a)\beta(1) + \alpha(2-a)\beta(2) + \cdots\\
 &+ \alpha(n-1)\beta(n-1+a)\\
 &= \sum_{j=0}^{n-1} \beta(j)\alpha(j-a),\\
 \sum_{i=0}^{n-1} \alpha(i)\beta(i-a) &= \alpha(0)\beta(-a) + \alpha(1)\beta(1-a) + \cdots + \alpha(n+a-1)\beta(n-1)\\
 &+ \alpha(a)\beta(0) + \alpha(1+a)\beta(1) + \alpha(2+a)\beta(2) + \cdots \\
 &+ \alpha(n-1)\beta(n-1-a)\\
 &= \sum_{j=0}^{n-1} \beta(j)\alpha(j+a).
\end{align*}
\normalsize
Thus
\begin{align*}
(\alpha^1\cup\beta^1)(f) &= - \sum_{a=1}^{\lfloor\frac{n-1}{2} \rfloor} \bigg( \frac{1}{2} - \frac{a}{n} \bigg) \sum_{i=0}^{n-1} \beta(i)\big(\alpha((i+a)\%n) - \alpha((i-a)\%n)\big)\\
&= -(\beta^1\cup\alpha^1)(f).
\end{align*}
Hence the polygonal cup product is skew--commutative on the cochain level.

Regarding the announced form of \textit{associativity} in item 4, the case of the product of three 0--forms corresponds to multiplying three real numbers associated to each vertex, which is indeed an associative operation. For the other degrees, we have learned in the proof of Theorem \ref{theo:PolygonalCupIsCup} that if $K$ is connected and a 0--form $\alpha^0$ is constant (which is equivalent to closed), i.e., $\alpha^0(v)= \rho\in\mathbb{R}\;\forall v\in K$, then the cup product of $\alpha^0$ with any $q$--form $\omega^q$ is equal to multiplying $\omega^q$ by the real number $\rho$. This fact yields:
 \[
   (\alpha^0\cup\beta^k)\cup\gamma^l = (\rho\beta^k)\cup\gamma^l = \rho(\beta^k\cup\gamma^l) =\alpha^0\cup(\beta^k\cup\gamma^l),
 \]
 which proves the statement of item 4.

We prove the item 5 for each type separately. Let $e=(v_0,v_1)$, then
\[
\begin{array}{rl}
\big(\alpha^0\cup(\beta^0\cup\gamma^1)-(\alpha^0\cup\beta^0)\cup\gamma^1\big)(e)
=& \frac{\alpha(v_0) + \alpha(v_1)}{2}\frac{\beta(v_0) + \beta(v_1)}{2}\gamma(e)\\
&- \frac{\alpha(v_0)\beta(v_0)  + \alpha(v_1)\beta(v_1)}{2}\gamma(e)\\
=& \frac{\gamma(e) }{4}(\alpha(v_1) - \alpha(v_0))
(\beta(v_0)  - \beta(v_1)).
\end{array}
\]
If $d\alpha$ is bounded, then $\alpha$ is Lipschitz continuous, thus
$|\alpha(v_1) - \alpha(v_0)|\to 0$ for $\norm{v_1-v_0}\to 0$, and similarly for
$\beta$.\\
Let $f=(v_0,\dots,v_{p-1})$ and $e_i = (v_i,v_{(i+1)\%p})$, then after some algebraic manipulation we get
\begin{multline*}
\big(\alpha^0\cup(\beta^1\cup\gamma^1)-(\alpha^0\cup\beta^1)\cup\gamma^1\big)(f)=\\
\sum_{i=0}^{p-1}\sum_{j=0}^{p-1}
\big(\alpha(v_j)-\alpha(v_i) + \alpha(v_j)-\alpha(v_{i+1})\big)
\frac{\beta(e_i)}{2p}\sum_{a=1}^{\lfloor\frac{p-1}{2} \rfloor} \bigg( \frac{1}{2} - \frac{a}{p} \bigg) \big(\gamma(e_{i+a}) - \gamma(e_{i-a})\big).
\end{multline*}
If $d\alpha$ is again bounded, then $\alpha$ is Lipschitz continuous, therefore
$|\big(\alpha(v_j)-\alpha(v_i) + \alpha(v_j)-\alpha(v_{i+1})\big)|\leq|\alpha(v_j) - \alpha(v_i)| + |\alpha(v_j)-\alpha(v_{i+1})\big)|\to 0$ for $\norm{v_j-v_i}\to 0$ and $\norm{v_j-v_{i+1}}\to 0$. That is, $\big(\alpha^0\cup(\beta^1\cup\gamma^1)-(\alpha^0\cup\beta^1)\cup\gamma^1\big)(f)$ converges to a constant 2--form, which is zero on all faces of the mesh.
\end{proof}

\section{Conclusion}\label{sec13}
We have presented a definition of a new discrete wedge product on polygonal 2--pseudomanifolds (Definition \ref{def:PolygonalCup}) and proved that this polygonal wedge product satisfies the Leibniz product rule with discrete exterior derivative (Lemma \ref{lemma:LeibnizPolygonalCup}). We have also shown that the polygonal wedge product is a cup product (Theorem \ref{theo:PolygonalCupIsCup}). Moreover, in Proposition \ref{prop:PolygonalCup} we have proved that the polygonal wedge product shares further properties with the wedge product of differential forms, such skew--commutativity, and that it is \textit{almost} associative.

\section*{Declaration of competing interest}
The author declares that she has no competing interests.

\section*{Acknowledgement}
The author was supported in part by the Czech Science Foundation Grant 25-15714S
and in part by the ERC Starting Grant ERC-2022-StG 101075632.

\bibliographystyle{elsarticle-harv}
\bibliography{PolyWedge}

\end{document}